\documentclass[12pt, twoside]{article}   	
\usepackage{geometry}                		
\usepackage{graphicx}				
\usepackage{amssymb}
\usepackage{authblk}
\usepackage{url}

\title{Sir Michael Atiyah, a Knight Mathematician\\
{ \small A tribute to Michael Atiyah, an inspiration and a friend\footnote{“This is Part 1 of a two-part series on Sir Michael Atiyah. Part 2 will be included in the December issue of Notices of the AMS.”}}}

\author{Alain Connes and Joseph Kouneiher}

\def\Ind{\mathop{\rm Ind}\nolimits}
\def\ch{\mathop{\rm ch}\nolimits}
\def\Diff{\mathop{\rm Diff}\nolimits}
\date{}							

\begin{document}
\maketitle

\noindent Sir Michael Atiyah was considered  one of the world's foremost mathematicians.  He is best known for his work in algebraic topology and the codevelopment of a branch of mathematics called topological K-theory together with  the Atiyah-Singer index theorem for which he received Fields Medal (1966).  He also received the Abel Prize (2004) along with Isadore M. Singer for their discovery and proof of the index theorem, bringing together topology, geometry and analysis, and for their outstanding role in building new bridges between mathematics and theoretical physics. Indeed,  his work has helped theoretical physicists to advance their understanding of quantum field theory and general relativity.

\vspace{0.2cm}

\noindent Michael's approach to mathematics was based primarily on the idea of finding new horizons and opening up new perspectives. Even if the idea was not validated by the mathematical criterion of  proof at the beginning, ``the idea  would become rigorous in due course, as happened in the past when Riemann used analytic continuation to justify Euler's brilliant theorems." For him an idea was justified by the new links between different problems which it illuminated. 
Our experience with him is that,  in the manner of an explorer,  he adapted to the landscape he encountered on the way  until he conceived a global vision of the setting of the problem.

\noindent Atiyah describes here \footnote{The quotations of Michael Atiyah presented in this paper have as sources personal exchanges with Michael or interviews that Michael gave on various occasions.} his way of doing mathematics\footnote{To learn more about Atiyah approach to mathematics, we refer to the general section of Atiyah's collected works \cite{atiyahcollectedpapers}.} :

\begin{quote}
Some people may sit back and say, I want to solve this problem and they
sit down and say, ``How do I solve this problem?" I don't. I just move around in
the mathematical waters, thinking about things, being curious, interested, talking to
people, stirring up ideas; things emerge and I follow them up. Or I see something
which connects up with something else I know about, and I try to put them together
and things develop. I have practically never started off with any idea of what I'm going to be doing or where it's going to go. I'm interested in mathematics; I talk, I
learn, I discuss and then interesting questions simply emerge. I have never started off with a particular goal, except the goal of understanding mathematics.
\end{quote}

\noindent We could describe Atiyah's  journey in mathematics by saying he  spent the first half of his career connecting mathematics to mathematics, and the second half connecting mathematics to physics.

\section{Building Bridges}

\noindent  Michael Atiyah's early education took place in English schools in Khartoum, Cairo and Alexandria, but his love and preference for geometry began as early as his two years at Manchester Grammar School which he spent preparing  to take the Cambridge scholarship examinations \cite{hitchin2}:

\begin{quote}
    I found that I had to work very hard to keep up with the class and the competition was stiff. We had an old-fashioned but inspiring teacher who had graduated from Oxford in 1912 and from him I acquired a love of projective geometry, with its elegant synthetic proofs, which has never left me. I became, and remained, primarily a geometer though that word has been reinterpreted in different ways at different levels. I was also introduced to Hamilton's work on quaternions, whose beauty fascinated me, and still does. (The Abel Prize 2004 lecture)
    \end{quote}

\noindent He won a scholarship to Trinity College, Cambridge in 1947. He read Hardy and Wright's Number Theory  and some articles on group theory during his two-year National Service. Thanks to his exceptional talent, he came out ranked first and wrote his first paper ``{\it A note on the tangents of a twisted cubic}" (1952) while still an undergraduate.

\noindent From 1952 to 1955 Atiyah undertook  three years of research at Trinity College, Cambridge, with William V. D. Hodge as supervisor.  He obtained his doctorate in 1955 with his thesis {\it ``Some Applications of Topological Methods in Algebraic Geometry"}. In fact, Atiyah hesitated between J.A.Todd and  Hodge:

\begin{quote}
When I came to decide on my graduate work I oscillated between Todd and Professor W. V. D. Hodge who represented a more modern approach based on differential geometry. Hodge's greater international standing swung the balance and, in 1952, I became his research student (\cite{atiyahcollectedpapers} vol.I). 
\end{quote}

\noindent Speaking of the work for his thesis, Atiyah said \cite{atiyahcollectedpapers}:

\begin{quote}
    I'd come up to Cambridge at a time when the emphasis in geometry was on classical projective algebraic geometry of the old-fashioned type, which I thoroughly enjoyed. I would have gone on working in that area except that Hodge represented a more modern point of view - differential geometry in relation to topology; I recognized that. It was a very important decision for me. I could have worked in more traditional things, but I think that it was a wise choice, and by working with him I got much more involved with modern ideas. He gave me good advice and at one stage we collaborated together. There was some recent work in France at the time on sheaf theory. I got interested in it, he got interested in it, and we worked together and wrote a joint paper which was part of my thesis. That was very beneficial for me.
\end{quote}

\noindent Atiyah became interested in analytic fiber bundles following his encounter with Newton Hawley who was visiting Cambridge. Reading the ``Comptes Rendus" allowed him to follow the development of sheaf theory in France. His single-author paper ``{\it Complex fibre bundles and ruled surfaces}" \cite{atiyah1955a} for which he received the Smith's Prize was motivated by the treatment of Riemann-Roch from the sheaf-theory\footnote{It seems that Serre's letter to Andr\'e Weil that gave the sheaf-theory treatment of Riemann-Roch for an algebraic curve, presented to him by Peter Hilton, motivated Atiyah's paper \cite{atiyah1955a}.} point of view and the classification of fiber bundles over curves : ``{\it it looked at the old results on ruled surfaces from the new point of view}".

\vspace{0.1cm}

\noindent Atiyah published two joint papers with  Hodge, {\it`` Formes de seconde esp\`ece sur une vari\'et\'e alg\'ebrique }" \cite{atiyah1954} and ``{\it Integrals of the second kind on an algebraic variety}" \cite{atiyah1954}. Those papers and \cite{atiyah1955a} granted him  the award of a Commonwealth Fellowship to visit  the Institute for Advanced Study in Princeton for the 1955-56 session.  On this occasion Atiyah met Jean-Pierre Serre, Friedrich Hirzebruch, Kunihiko Kodaira, Donald Spencer, Raoul Bott,  Isadore Singer, and others. During his stay at the IAS Atiyah attended Serre's seminars, which   were a source of inspiration and motivated  his following papers,  in which he showed how the study of vector bundles on spaces could be regarded as the study of a cohomology theory called $K$-theory.

\subsubsection*{$K$-theory}

The historical origin of $K$-theory\footnote{ It takes its name from the German {\it Klasse}, meaning "class".  The work of J. H. C. Whitehead or Whitehead torsion can be considered ultimately as the other historical origin.} can be found in Alexander Grothendieck's work on coherent sheaves on an algebraic variety $X$ \cite{grothendieck1957}. Grothendieck used it to formulate his Grothendieck-Riemann-Roch theorem. The idea is to define a group using isomorphism classes of sheaves as generators of the group, subject to a relation that identifies any extension of two sheaves with their sum. The resulting group is called $K(X)$ when only locally free sheaves are used, or $G(X)$ when all are coherent sheaves\footnote{$K(X)$ has cohomological behavior and $G(X)$ has homological behavior.}.  

\begin{quote}
Having been exposed, at the {\it Arbeitstagung}, to many hours of Grothendieck expounding his generalization of the Grothendieck-Riemann-Roch theorem, I was playing about with his formulae for complex projective space. From Ioan James (at that time a colleague in Cambridge) I had heard about the Bott periodicity theorems and also about stunted projective spaces. I soon realized that Grothendieck's formulae led to rather strong results for the James problems. Moreover the Bott periodicity theorems fitted in with the Grothendieck formalism, so that one could draw genuine topological conclusions. It was this which convinced me that a topological version of Grothendieck's $K$-theory, based on the Bott periodicity theorem, would be a powerful tool in algebraic topology  (\cite{atiyahcollectedpapers} vol2).
\end{quote}

\noindent Using the same construction applied to vector bundles in topology,  Atiyah and Hirzebruch defined $K(X)$ for a topological space\footnote{Atiyah's joint paper with Todd was the starting point of this work \cite{atiyah1960}.} $X$ in 1959 \cite{atiyahhirzebrucha}, and using the Bott periodicity theorem they made it the basis of an extraordinary cohomology theory \cite{atiyahhirzebruchd}. It played a major role in the second proof of the Index Theorem\footnote{This approach led to a noncommutative $K$-theory for ${\mathbb C}^*$-algebras. 
}.

\noindent The Atiyah--Hirzebruch spectral sequence relates the ordinary cohomology of a space to its generalized cohomology theory. Atiyah applied the approach  to finite groups $G$ and showed \cite{{atiyah1961a}} that for a finite group $G$, the $K$-theory of its classifying space, $BG$, is isomorphic to the completion of its character ring :

    \begin{equation}
    K(BG)\cong R(G)^{{\wedge }}.
    \label{K(BG)}
\end{equation}

\noindent The same year  Atiyah and Hirzebruch established the relation between the representation ring $R(G)$ of a compact connected Lie group and $K$-theory of its classifying space \cite{atiyahhirzebruchd}.

\noindent Later on the result could be extended to all compact Lie groups using   results from Graeme Segal's thesis. However a simpler and more general proof was produced by introducing equivariant $K$-theory, i.e. equivalence classes of $G$-vector bundles over a compact $G$-space $X$ \cite{atiyahsegal}. In this paper Atiyah and Segal showed that under suitable conditions the completion of the equivariant $K$-theory of $X$ is isomorphic to the ordinary $K$-theory of a space, $X_{G}$, which fibers over $BG$ with fibre $X$:

$$  K_{G}(X)^{{\wedge }}\cong K(X_{G}).$$

\noindent The original result (\ref{K(BG)}) can be recovered as a corollary by taking $X$ to be a point. In this case the left hand side reduced to the completion of $R(G)$ and the right to $K(BG)$.

\noindent In \cite{atiyah1962} Atiyah gave a $K$-theory interpretation of the {\it symbol}\footnote{The principal symbol of a linear differential operator $\Sigma_{ \mid\alpha\mid \leq m} a_{\alpha} (x) \partial_x^{\alpha}$ is by definition the function $\Sigma_{ \mid\alpha\mid \leq m} a_{\alpha} (x) (i\xi)^{\alpha}$  obtained from the differential operator (of a polynomial) by replacing each partial derivative by a new variable.  At this point, the vector $\xi = (\xi_1, \ldots, \xi_n)$ is merely a formal variable. 

An essential and interested fact is that if one interprets $(x, \xi)$ as variables in the cotangent bundle then the principal symbol is an invariantly defined function on $ T^*X$, where $X$ is the manifold on which the operator is initially defined, which is homogeneous of degree $m$ in the cotangent variables. the symbol captures locally some very strong properties of the operator. For example, an operator $𝐿$ is called elliptic if and only if the symbol is invertible (whenever $\xi \ne 0$).} of an elliptic operator. This work induced his interest in elliptic operators which then became a dominant theme in his work, so that his later papers on $K$-theory appeared interspersed chronologically with his papers on index theory\footnote{Atiyah and Singer's interest in the elliptic operator questions was dependent on the fact that  real $K$-theory behaves differently from complex $K$-theory, and requires a  different notion of reality. The essential point is that the Fourier transform of a real-valued function satisfies the relation $f(-x) = \bar{f(x)}$. Their new version of real $K$-theory involved spaces with involution.}. Using $K$-theory Atiyah was able to give a simple and short proof of the Hopf invariant problem. This involved the $\psi^k$ operations which were introduced by Adams in his subsequent work on the vector field problem.

\noindent Algebraic $K$-theory was successfully developed  by Daniel Quillen using homotopy theory in 1969 and 1972, following a trend of ideas which can be traced back to the work of Whitehead in 1939 \cite{quillen2}. A further extension was developed by Friedhelm Waldhausen in order to study the algebraic $K$-theory of spaces, which is related to the study of pseudo-isotopies. Much modern research on higher $K$-theory is related to algebraic geometry and the study of motivic cohomology. The corresponding constructions involving an auxiliary quadratic form received the general name $L$-theory,   a major tool of surgery theory. In string theory, the $K$-theory classification of Ramond--Ramond field strengths and the charges of stable $D$-branes was first proposed in 1997.

\subsubsection*{The index theorem}

The concept of a linear operator, which together with the concept of a vector space is fundamental in linear algebra, plays a role in very diverse branches of mathematics and physics, above all in analysis and its applications. 
Up to the beginning of the 20th century the only linear operators that had been systematically studied were those between finite-dimensional spaces over the fields ${\mathbb R}$ and ${\mathbb C}$. For instance, the classical theory of Fredholm integral operators goes back at least to the early 1900s \cite{fredholm}.  Fredholm essentially proved (and D. Hilbert perceived this "geometrical" background) that a linear operator of the form $ {\mathbb I} + {\mathbb A}$, with compact ${\mathbb A}$, has index [$dim\,\, ker({\mathbb A}) - dim \,\,coker({\mathbb A})$] equal zero.

\noindent Special forms of index formulas were known even earlier, for example, the Gauss-Bonnet theorem and its multi-dimensional variants.  Indeed, the precursors of the index theorem were formulas for Euler characteristics  -- alternating sums of dimensions of cohomology spaces. Putting the even ones on one   side and the odd ones on the other and using Hodge theory this could be seen in more general terms as the difference between the dimension of the kernel and the cokernel of an elliptic operator\footnote{\noindent  In general,  the dimension of the cokernel and kernel of an elliptic operator are extremely hard to evaluate individually -- they even vary under deformation -- but  the index theorem shows that we can usually at least evaluate their difference. In applications, it is sometimes possible to show that the cokernel is zero and then we have  a formula for the dimension of the null space without solving the equation explicitly. }. Subsequently a number of generalizations of index formulas were obtained for objects of a more complex nature.

\noindent The first "infinite-dimensional" observations, concerned also with general fields, were made by O. Toeplitz \cite{toeplitz}. As a rule, linear operators between infinite-dimensional spaces are studied under the assumption that they are continuous with respect to certain topologies. Continuous linear operators that act in various classes of topological vector spaces, in the first place Banach and Hilbert spaces, are the main object of study of linear functional analysis

\noindent Historically, the work of Fritz Noether in analysis in 1921 \cite{noether1} exhibit already a case where the index is different from zero. Noether gave a formula for the index in terms of a winding number constructed from data defining the operator. This result was generalized by G. Hellwig, I.N. Vekua and others (see \cite{vekua}),

\noindent  Later Israel Gel'fand noticed in 1960 \cite{gelfand1960} that the {\it index}, the difference in these dimensions is invariant under  homotopy. This suggested that the index of a general elliptic operator should be expressible in terms of topological data, and Gel'fand asked for a formula for it by means of topological invariants. Gel'fand proposed that the index of an elliptic differential operator (with suitable boundary conditions in the  presence  of  a  boundary)  should  be  expressible  in  terms  of  the  coefficients  of highest order part (i.e., the principal symbol) of the operator, since the lower order parts provide only compact perturbations which do not change the index.  Indeed, a  continuous,  ellipticity-preserving  deformation  of  the  symbol  should  not  affect the index\footnote{In some sense index theory is about regularization, more precisely, the index quantifies the defect of an equation, an operator, or a geometric configuration from being regular.  Index  theory  is  also  about perturbation  invariance,  i.e.,  the  index  is  a meaningful quantity stable under certain deformations and apt to store certain topological or geometric information. Most important for many mathematicians, the index interlinks  quite  diverse mathematical fields, each with its own very distinct research tradition.}, and so Gel'fand noted that the index should depend only on a suitably defined homotopy class of the principal symbol.  The hope was that the index of an elliptic operator could be computed by means of a formula involving only the topology of the underlying domain (the manifold), the bundles involved, and the symbol  of  the  operator.

\noindent In  early  1962,  M.F.  Atiyah  and  I.M.  Singer  discovered the (elliptic) Dirac operator in the context of Riemannian geometry and they were concerned by the proof  of the fact that  the $\hat{A}$-genus  of  a  spin  manifold  is  the  index of  this  Dirac  operator.  They were concerned also by the  geometric meaning of the $\hat{A}$-genus. Hirzebruch, with  Borel, had proved that for  a spin manifold the index  was an integer. 
  
 \begin{quote}
 [The index theory of elliptic operators] had its origins in my work on $K$-theory with Hirzebruch and the attempt to extend the  Hirzebruch-Riemann-Roch theorem into differential geometry  \cite{atiyahhirzebrucha, atiyahhirzebruchd}. We had already shown that the integrality of the Todd genus of an almost complex manifold and the $\hat{A}$-genus of a spin-manifold could be elegantly explained in terms of $K$-theory. For an algebraic variety the Hirzebruch-Riemann-Roch theorem went one step further and identified the Todd genus with the arithmetic genus or Euler characteristic of the sheaf cohomology. Also the $L$-genus of a differential manifold, as proved by Hirzebruch, gave the signature of the quadratic form of middle dimensional cohomology and, by Hodge theory, this was the difference between the dimensions of the relevant spaces of harmonic forms. As Hirzebruch himself had realized it was natural therefore to look for a similar analytical interpretation of the $\hat{A}$-genus. The cohomological formula and the associated character formula clearly indicated that one should use the spin representations  (\cite{atiyahcollectedpapers} vol4). 
 \end{quote}

\noindent  It was Stephen Smale who turned their attention to Gel'fand's general program,

\begin{quote} 
 Once we had grasped the significance of spinors and the Dirac equation it became evident that the $\hat{A}$-genus had to be the difference of the dimensions of positive and negative harmonic spinors. Proving this became our objective. By good fortune Smale passed through Oxford at this time and, when we explained our ideas to him, he drew our attention to a paper of Gel'fand on the general problem of computing the index of elliptic operators (\cite{atiyahcollectedpapers} vol4).
\end{quote} 

\noindent The generalization in 1962 of F. Hirzebruch's proof of the Hirzebruch-Riemann-Roch theorem \cite{hirzebruch} and the proof of the desired index formula by Atiyah and Singer were based in fact on work by Agranovic, Dynin, and others analysts\footnote{\noindent By investigating the particular case of Dirac operator Atiyah and Singer had an advantage over the analysts investigating  the index problem, and mostly because they knew already what had to be the answer namely the $\hat{A}$-genus. In fact this case encompasses virtually all the global topological applications by twisting with an appropriate vector bundle.},  particularly  those  involving pseudo-differential operators.

\noindent In their first published proof\footnote{The proof in \cite{atiyah1968a} is based on the idea of the invariance of the index under homotopy which implies that the index (say, the analytic index) of an elliptic operator is stable under continuous changes of its principal symbol while maintaining ellipticity. Using this fact, one finds that the analytical index of an elliptic operator transforms predictably under various global operations such as embedding and extension. Using K-theory and Bott periodicity, we can construct, from the symbol of the elliptic operator, a topological invariant (the topological index) with the same transformation properties under these global operations is constructed.} \cite{atiyah1968a} they used $K$-theory instead of cohomology and replaced the cobordism\footnote{The details of this original proof involving cobordism actually first appeared in \cite{palais}.} theory of the first proof with arguments akin to those of Grothendieck. The proof  is based on embedding a manifold in Euclidean space and then transferring the problem to one on  Euclidean space by a suitable direct image construction. It worked purely in a $K$-theory context and avoided rational cohomology. They used this then to give proofs of various generalizations in their papers  \cite{atiyah1968a}.

\noindent We can then verify that a general index function having these properties is unique, subject to normalization. To deduce the Atiyah-Singer Index theorem (i.e., analytic index = topological index), it then suffices to check that the two indices are the same in the trivial case where the base manifold is just a single point. In other words, the index problem concerns the computation of the (analytical) index of $D$ using only its highest order term, the  symbol. The symbol defines a homomorphism between vector bundles and naturally defines a class in K-theory.  Topological data derived from the symbol and the underlying manifold defines a rational number, the topological index, and we have the theorem \cite{ atiyah1963, palais}

  \begin{quote}  
  {\it The analytical index of $D$ is equal to its topological index}.
  \end{quote}

\noindent The key point is the fact that we can usually evaluate the topological index explicitly,  so this makes it possible to evaluate the analytical index. 
 
 \vspace{0.2cm}

\noindent Atiyah, Bott, and Patodi gave a second proof of the index theorem using the heat equation \cite{atiyah-bott-patodi}. The heat kernel method had its origins in the work of Mckean in the late 1960s \cite{mckean} drawing on \cite{mipl} and was pioneered in the works of Patodi \cite{patodi} and Gilkey \cite{gilkey}. The method can also largely applies for more general elliptic pseudo-differential operators. 

\noindent  In the heat kernel approach, the invariance of the index under changes in the geometry of the operator is a consequence of the index formula itself more than a means of proof. Moreover, it can be argued that in some respects the $K$-theoretical embedding/cobordism methods are more subtle and direct. For instance, in some circumstances the index theorem involves torsion elements in $K$-theory that are not detectable by cohomological means, and hence are not computable in terms of local densities produced by heat asymptotics \cite{lami}.

\noindent More importantly, the heat kernel approach exhibits the index as just one of a whole sequence of spectral invariants appearing as coefficients of terms of the asymptotic expansion of the trace of the relevant heat kernel.

\noindent  The  third proof, in  the  early  seventies,  aimed  for  a  more   direct  transition  between  the  analysis  and  characteristic  classes  by  expressing these as curvature integrals. This local proof of the index theorem  relies on a detailed asymptotic analysis of the heat equation associated to a Dirac operator.

\noindent This is analogous to certain intrinsic proofs of the Gauss-Bonnet theorem, but the argument doesn't provide the same kind of intuition that the $K$-theory argument does. The basic strategy of the local argument is to invent a symbolic calculus for the Dirac operator which reduces the theorem to a computation with a specific example. 

\noindent Many subsequent developments and applications have appeared since then in the work of Widom \cite{wid}, Getzler \cite{get, begeve} and Pflaum \cite{pf}. Getzler Draws on a version of the quantum-mechanical harmonic oscillator operator and a coordinate calculation directly to  produce the $\hat{A}$-genus (the appropriate ``right-hand side" of the index theorem for the Dirac operator)

\noindent The more  recent  proofs  of  the index  theorem which  are  motivated  by  supersymmetry  provide  a  rationale  for  the  role  which   functions  like ${x}/{(1-e^{-x})}$ and ${x}/{\tanh x}$ play  in  sorting out the combinations of  characteristic  classes  which  occur  in  the  index  formula.

\noindent In the first years the index  theorem  was  applied in  the  area  of  topology  or  number  theory.  A  typical  number  theoretical  application  would  involve  taking  a  group  action  on  a  manifold  that  was  so  well  analyzed  that  both  the  analytical  side  and  the  topological  side  could  be   calculated   independently.   The   index   theorem   then   produced   an   identity, sometimes well known, at other times new. 

\noindent Today the index theorem plays an essential role in partial  differential  equations, stochastic processes, Riemannian  geometry, algebraic geometry,  algebraic topology and mathematical physics.

 \vspace{0.2cm}

\noindent Since the appearance of the Atiyah-Singer index theorem  generalizations have been obtained encompassing many new geometric situations. Among them some constructions involve non-Fredholm operators and use in a decisive way an extra ingredient consisting of von Neumann algebras and the Murray-von Neumann dimension theory. While the numerical index of elliptic operators on non compact manifolds does not make much sense, some useful assumptions can enable us to associate them to appropriate type $II$ von Neumann algebras and to give a sense to the Murray-von Neumann's index, particularly in the case of a manifold with an infinite fundamental group. The use of von Neumann algebras was quite natural and led to concrete nontrivial results.

\begin{quote}
Since I was not an expert on von Neumann algebras I attempted in my presentation (at the meeting in honour of Henri Cartan) to give a simple, elementary and essentially self-contained treatment of the results. Later on in the hands of Alain Connes, the world expert on the subject, these simple ideas were enormously extended and developed into a whole theory of linear analysis for foliations (\cite{atiyahcollectedpapers} vol4 and \cite{connes1985b}).
\end{quote}

\noindent The index theorem has had many variants, for example, to  the extension of Hodge theory on combinatorial and Lipschitz manifolds and Atiyah-Singer's signature operator to Lipschitz manifolds.
For any quasi-conformal manifold there exists a local construction of the Hirzebruch-Thom characteristic classes \cite{alain1}. This work \cite{alain1} provides local constructions for characteristic classes based on higher dimensional relatives of the measurable Riemann mapping theorem in dimension two and
the Yang- Mills theory in dimension four. At the same time, they provide also an effective construction of the rational Pontrjagin classes on topological manifolds. In 1985 Teleman provided in his paper a link between index theory and Thom's original construction of the rational Pontrjagin classes.

\subsubsection*{From the Atiyah-Singer index theorem to noncommutative geometry}

 \noindent   Besides topological $K$-theory, Michael Atiyah made two discoveries that played a major role in the elaboration of the fundamental  paradigm of noncommutative geometry (NCG). This theory is an extension of Riemannian geometry to spaces whose coordinate algebra ${\mathcal{A}}$ is no longer commutative. A geometric space is encoded spectrally by  a representation of ${\mathcal{A}}$ as operators in a Hilbert space ${\mathcal{H}}$ together with an unbounded self-adjoint operator $D$. 
Consider an $n$-dimensional compact manifold $M$ with Riemannian metric, and  spin structure. This means that there is a Hilbert space ${\mathcal{H}}$  of {\it spinor fields}  $\psi(x)$ on $M$ (i.e. ``sections of the spinor bundle"), on which both the algebra of functions on $M$\footnote{A function $f(x)$ acts by multiplication $(f.\psi)(x):=f(x)\psi(x)$ on spinors.} and the Dirac operator\footnote{In local coordinates $D = i\gamma_{\mu}\partial x^{\mu}$,   $\gamma_{\mu}$ the Dirac $\gamma$-matrices, which satisfy the wellknown condition  $\gamma_{\mu}\gamma_{\nu} + \gamma_{\nu}\gamma_{\mu}= 2g_{\mu\nu}$.} $D$ act. The ``spectral triple" $({\mathcal{A}},{\mathcal{H}},D)$ faithfully encodes the original Riemannian geometry and opens the door to the wide landscape of new geometries that cover leaf spaces of foliations, the fine structure of space-time as well as the infinite-dimensional geometries of quantum field theory.

\noindent Behind this concept of spectral triple is the realization by Atiyah \cite{atiyah1970} that the cycles in the homology theory, called $K$-homology, which is dual to $K$-theory are naturally given by a Fredholm representation of the algebra $A=C(X)$ of continuous functions on the compact space $X$. Moreover, one owes to Atiyah and Singer the identification of the $KO$-homology class of the Dirac operator of a spin manifold as the fundamental cycle realizing  Poincar\'e duality in $KO$-homology. This Poincar\'e duality had been shown by Dennis Sullivan to be a key property of manifolds, much deeper than Poincar\'e duality in ordinary homology since for instance this fundamental cycle determines the Pontryagin classes of the manifold. 
One consequence of the duality is that it reduces the index problem to the computation of the index of twisted Dirac operators\footnote{Notice that a reduction is already worked out in \cite{atiyah-bott-patodi}. The idea basically is that the classical operators are sufficiently numerous to generate, in a certain topological sense, all elliptic operators. }. In this case and for any hermitian vector bundle $E$ over an even-dimensional   compact spin manifold $M$ the index formula is $\Ind(\partial_E^+) = < \ch(E) {\hat{A}}(M), \, [M]>$
 \noindent  where $\partial_E^+$ is the positive part of the self-adjoint Dirac operator twisted by the
bundle $E$, $\ch(E)$ is the Chern character of $E$ and $\hat{A}(M)$ is the $\hat{A}$-polynomial of the
tangent bundle $TM$ of $M$ \cite{benameur}. 
 Treating this index
as a pairing between $E$ and the operator $\partial$ \cite{atiyah1970},  defines the Chern character of the Dirac operator as the homology class
\begin{equation}
\ch(\partial) =  {\hat{A}}(M) \star [M]\, \label{3}
\end{equation}

 \noindent  From equation [\ref{3}] we can  interpret the index theorem as a homological pairing 
 \begin{equation}
\Ind(\partial_E^+) = < \ch(E), \ch(\partial)>\,.\label{4}
\end{equation}
 The key fact is that this pairing between $K$-theory and $K$-homology makes sense in the general noncommutative case and can be computed using a natural map from $K$-homology to cyclic cohomology as well as a ``Chern character" from $K$-theory to cyclic homology. Then pairing these two elements from cyclic cohomology and homology gives the same value as if we start from $K$-theory and $K$-homology.

\noindent Equation [\ref{4}] represents the noncommutative geometry point of view of the index theorem. To make the computation effective, the idea is to obtain a Local Index formula. Given a spectral triple  the Connes-Moscovici local index theorem\footnote{The term local refers to the fact that in the classical situation the index formula would involve integration and could be computed from  local data. Another point of view is that ``we are working in momentum space and local means therefore asymptotic".} is then a complete solution of this general index problem. The first highly nontrivial spectral triple is the example which encodes the $\Diff$-equivariant index theory. The identification of the local terms appearing in the index theorem for such a triple requires the definition of  an appropriate cyclic homology for Hopf algebras \cite{connespreprint} encoding the  characteristic numbers.

\section{A new era in mathematical physics}

 \noindent  The ``recent " interaction between physics and mathematics, and particularly geometry, is mainly due to
the fact that physicists started exploring complicated nonlinear
models of a geometrical character as possible explanations for fundamental
physical processes \cite{kouneiher1}. Differential geometry had been
closely associated with physics ever since the introduction of
Einstein's theory of general relativity but quantum theory, the essential
ingredient for the  physics of elementary particles, had not impacted the type of  geometry involving global topological features. One of the first occurrence of topology\footnote{One of the earliest occurrences of some topological noteworthiness can be found in Euler's 1736 solution of the bridges of K\"{o}nigsberg problem. In solving one of the first problems in combinatorial topology, Euler gave birth to what we now call graph theory. A close connection between physics and  graph theory occurs in the work of Kirchhoff in 1847 on his two famous laws for electric circuits. 

\noindent In the nineteenth century we encounter more substantial examples of physical phenomena with topological content. In Gauss's 1833 work on electromagnetic theory, found in his Nachlass (Estate), the topology concerned the linking number of two curves. Another work was in the theory of vortex atoms proposed by Lord Kelvin (alias W. H. Thomson) in 1867. Kelvin was influenced by an earlier fundamental paper by Helmholtz (1858) on vortices, and a long and seminal paper of Riemann (1857) on Abelian functions (see \cite{nash}).

 \noindent Maxwell on the state of the art in the nineteenth century:

 \begin{quote}

 {\it It was the discovery by Gauss of this very integral, expressing the work done on a magnetic pole while describing a closed curve in presence of a closed electric current, and indicating the geometrical connexion between the two closed curves, that led him to lament the small progress made in the Geometry of Position since the time of Leibnitz, Euler and Vandermonde. We have now some progress to report, chiefly due to Riemann, Helmholtz, and Listing}.
 \end{quote}

 \noindent Poincar\'e's work on the {\it Analysis Situs} enlivened considerably the topological matters and its link to Physics.} in physics was Dirac's famous argument for the quantization of electric charge, and what we are now witnessing is essentially the non-abelian outgrowth of Dirac's initial idea. 
 \noindent  This intimate relation between geometry and physics\footnote{For more on the relation between Physics and Topology see \cite{nash}.} has had beneficial
consequences in both directions. Physicists have been able to
adopt and use sophisticated mathematical ideas and techniques without
which the elaboration of the physical theories would be greatly impeded.
Conversely mathematicians have used the insights derived from a
physical interpretation to break new ground in geometry.

\noindent Atiyah has been influential in stressing the role of topology in quantum field theory and in bringing the work of theoretical physicists to the attention of the mathematical community.

\vspace{0.2cm}

\noindent  In his paper with Bott \cite{atiyahbott}, they attempt to understand the results of V. K. Patodi on the heat equation approach to the index theorem, a topic which will be developed later in connection with theoretical physics. In their work on the Lefschetz formula for elliptic complexes they had described  the Zeta-function approach to the index theorem.  Singer and Mckean  in their paper \cite{mckean} took a step further, in the heat equation version, by concentrating on the Riemannian geometry. They suggest the possibility of cancellations leading directly to the Gauss-Bonnet from the Euler Characteristic\footnote{This was proved later by Patodi, who extend the result to deal with the Riemann-Roch theorem on K\"{a}hler manifolds.  On his side Gilkey  produced an alternative indirect approach depending on a simple characterization of the Pointryagin forms of a Riemannian manifold \cite{gilkey}. Later on, Atiyah, Bott and Patodi realized that Gilkey's results was a consequence of the Bianchi identities.}.

\noindent  Another topic that would have an impact on theoretical physics, concerned Atiyah's work on  {\it spectral asymmetry}. Aspects of spectral asymmetry appear in differential geometry, topology, analysis and number theory. These applications became  popular, especially after Witten's work on global anomalies brought the $\eta$-invariant into prominence. The $\eta$-invariant introduced in a joint paper with Patodi and Singer, was motivated, in part, by Hirzebruch's work on cusp singularities and in particular his result expressing the signature defect in terms of values of $L$-functions of real quadratic fields.
\vspace{0.2cm}

\noindent  In the 1970s and 1980s, methods gleaned from the index theorem unexpectedly played a role in the development of theoretical physics. One of the earliest and most important applications is 't Hooft's resolution of the $U(1)$ problem \cite{hooft}. This refers to the lack of a ninth pseudo-Goldstone boson (like the pions and Kaons) in QCD that one would naively expect from chiral symmetry breaking. There are two parts to the resolution. The first is the fact that the chiral $U(1)$ is anomalous. The second is the realization that there are configurations of finite action (instantons) which contribute to correlation functions involving the divergence of the $U(1)$ axial current. The analysis relies heavily on the index theorem for the Dirac operator coupled to the $SU(3)$ gauge field of QCD, the main remaining law in the theory of the strong interactions.  Albert Schwarz showed that some of the ingredients in the solution were best understood in terms of the Atiyah-Singer index theorem.There were also important applications in the 1990s to S-duality of $N =4$ Super Yang-Mills which involve the index theorem for the Dirac operator on monopole moduli spaces.

Describing all those activities and the first uses of the index theorem in physics Atiyah said

\begin{quote}
From 1977 onwards my interests moved in the direction of gauge theories and the interaction between geometry and physics. I had for many years had a mild interest in theoretical physics, stimulated on many occasions by lengthy discussions with George Mackey. However, the stimulus in 1977 came from two other sources. On the one hand Singer told me about the Yang-Mills equations, which through the influence of Yang were just beginning to percolate into mathematical circles. During his stay in Oxford in early 1977 Singer, Hitchin and I took a serious look at the self-duality equations. We found that a simple application of the index theorem gave the formula for the number of instanton parameters, $\ldots$ At about the same time A. S. Schwarz in the Soviet Union had independently made the same discovery. From this period the index theorem began to become increasingly familiar to theoretical physicists, with far-reaching consequences.

The second stimulus from theoretical physics came from the presence in Oxford of Roger Penrose and his group.$\ldots$ at that time I knew nothing about twistors, but on Penrose's arrival in Oxford we had lengthy discussions which were mutually educational. The geometry of twistors was of course easy for me to understand, since it was the old Klein correspondence for lines in $P_3$ on which I had been brought up. The physical motivation and interpretation I had to learn.
In those days Penrose made great use of complex multiple integrals and residues, but he was searching for something more natural and I realized that sheaf cohomology groups provided the answer. He was quickly converted and this made subsequent dialogue that much easier, since we now had a common framework (\cite{atiyahcollectedpapers} vol5).
\end{quote}

 \noindent  The theoretical physicists  working  with  gauge  theories   were  interested  in  calculating  the  number  of  zero  modes  of  a  differential  operator  in  terms  of  the  topological  charge  of  the  gauge  field.  Translating  the  language,  they  wanted  the dimension of $ker D$  in  terms  of  characteristic  classes,  and  this  is  precisely   what   the  index   theorem,  coupled  with  what  mathematicians  would  call  a  vanishing  theorem  for  $coker D$, provides. Naturally,  the mathematical  physicist  is  more  familiar  with  eigenvalues  of  operators  and  the  heat  kernel  than  with  characteristic  classes,  but  once  the  realization  came  that  each  side  had  something  to learn  from  the other  the most recent phase  of  index  theory was  formed.

 \subsubsection*{A man of intuition}

\noindent The term ``gauge theory" refers to a quite general class of quantum field theories used for the description of elementary particles and their interactions. These theories involve deep and interesting nonlinear differential equations, in particular, the Yang-Mills equations which have turned out to be particularly fruitful for mathematicians. 

\begin{quote}
Meanwhile I had been engrossed with the Yang-Mills equations in dimension 4. I realized that these questions were essentially trivial in dimension 2, but one day, walking across the University Parks with Bott it occurred to me that one might nevertheless be able to use the Yang-Mills equations to study the moduli spaces. The essential point was the theorem of Narasimhan and Seshadri stating that stable bundles arose from unitary representations of the fundamental group (\cite{atiyahcollectedpapers} vol5).

\end{quote}

Atiyah initiated much of the early work in this field, culminated by Simon Donaldson's work in four-dimensional geometry. Indeed, Andreas Floer and  Donaldson's works  provided  new insights into the topology of lower dimensional manifolds.  In three and four dimensions this was Floer homology and the Donaldson polynomials based on Yang-Mills theory\footnote{One can transpose the same ideas to maps from a surface into a symplectic manifold and express the Floer homology of $X$ and the symplectic homology of $X / \Sigma$, for a special class of Hamiltonians, in terms of absolute and relative Gromov-Witten invariants of the pair $(X, \Sigma)$ and some additional Morse-theoretic information. The key
point of the argument is a relation between solutions of Floer's equation and pseudoholomorphic
curves, both dfiened on the symplectization of a pre-quantization bundle over $\Sigma$.}.

\noindent Donaldson's work showed that, in contrast to other dimensions, there were some very subtle invariants in four dimensions. These invariants were preserved under smooth deformation but not under general continuous deformation \cite{donaldson1990}. 
Donaldson found his invariants by studying special solutions of the Yang--Mills equations the ``instantons" \cite{hitchin4}. These solutions are essentially localized to a small region in space--time, thereby describing an approximately instantaneous process.

\noindent Donaldson's work involved  ideas from nonlinear analysis of partial differential equations,
differential and algebraic geometry, and topology and a little less physics. Nevertheless, beside the influence of Atiyah-Hitchin-Singer's paper the the whole idea of studying a moduli space in this context (a space of connections up to gauge equivalence) had a physics origin in some sense.

\noindent Witten later showed that Donaldson's invariants could be
interpreted in terms of a quantum field theory and that this would have profound
consequences. Moreover, this field theory was a close cousin of the standard
theories used by particle physicists, except that it had a ``twist" that produced
topological invariants, not dependent on the intricate details of the underlying
geometry of space--time.  Beside Seiberg's previous work on dualities, this interpretation of Witten motivated Seiberg and Witten's work. They showed that the corresponding physical theory could be
solved in terms of a much simpler structure. Yang--Mills theory is fundamentally
based on a choice of non-Abelian Lie group, usually taken to be the group $SU(2)$.
The non-commutativity of this symmetry group leads to the nonlinearities of the
associated partial differential equations. However, in physics it was known that in
quantum theories these non-Abelian symmetries often manifest themselves only
at very short distance scales. At large distances the symmetry can be broken
to a much simpler Abelian group. For example, in the case of $SU(2)$ only the
circle group $U(1)$ of electromagnetism would appear together with possibly some
charged matter particles \cite {atiyah2010} (\cite{atiyahcollectedpapers} vol7). As Witten writes:

\begin{quote}

In the spring of 1987 Atiyah visited the Institute for Advanced Study and was
more excited than I could remember. What he was so excited about was Floer theory which he felt should be interpreted as the Hamiltonian formulation of a quantum field theory. Atiyah hoped that a quantum field theory with Donaldson polynomials as the correlation functions and Floer groups as the Hilbert spaces could somehow be
constructed by physics methods. 
$\ldots$ I was skeptical and though the idea was intriguing, I did not
pursue it until I was reminded of the question during another visit by Atiyah to the
Institute at the end of 1987. This time I dropped some of my prejudices and had
the good luck to notice that a simple twisting of $N=2$ supersymmetric Yang-Mills theory would give a theory with the properties that Atiyah had wanted \cite{witten1999}.
\end{quote}

\noindent Another topic that was developed following Atiyah's suggestions and indications was the connection with the  knot polynomial of Vaughan Jones.

\begin{quote}
The other problem that Atiyah recommended for physicists in the years was to understand the Jones knot polynomial via quantum field theory. It was from
him that I first heard of the Jones polynomial. There followed other clues about the Jones polynomial and physics \cite{witten1999}.
\end{quote}

``In fact, the relation between knot invariants and particles goes to the very
beginning of relativistic quantum field theory as developed by Feynman and
others in the 1940s. The basic idea is that, if we think of a classical particle
moving in space--time, it will move in the direction of increasing time. However,
within quantum theory the rules are more flexible. Now a particle is allowed to
travel back in time. Such a particle going backwards in time can be interpreted
as an anti-particle moving forwards in time. Once it is allowed to turn around,
the trajectory of a particle can form, as it were, a complicated knot in space--
time. The rules of quantum theory will associate to each such trajectory a
probability amplitude that describes the likelihood of this process actually
taking place.

In order to establish contact between this formulation of particle physics and knot theory, Witten had
to replace the usual four dimensions of space--time with three dimensions
and work with a special type of gauge theory based on the Chern--Simons
topological invariant. The full theory allowed for many choices, such as the gauge
group, the representation of the particles associated to the knot and coupling
constants.

Relating quantum field theory to knot invariants along these lines had
many advantages. First, it was fitted into a general framework familiar to
physicists. Second, it was not restricted to knots in three-dimensional Euclidean
space and could be defined on general three-dimensional manifolds. One
could even dispense with the knot and get an invariant of a closed three dimensional
manifold. Such topological invariants are now called quantum
invariants and have been extensively studied "  \cite {atiyah2010} (\cite{atiyahcollectedpapers} vol7).

\noindent  Atiyah related all of this to Witten's work on quantum mechanics and Morse theory, arguing that there should be a quantum field theory version of the story.  This was the birth of topological quantum field theory (TQFT), with Witten to work out a year later what the quantum field theory conjectured by Atiyah actually was. Soon after the first topological TQFT, Witten was to find a quantum theory with observables the Jones polynomial, now often called the Chern-Simons-Witten theory.

\noindent The theories of superspace and supergravity and the string theory of fundamental particles, which involves the theory of Riemann surfaces in novel and unexpected ways, were all areas of theoretical physics which developed using the ideas which Atiyah was introducing. 

\vspace{0.2cm}

\noindent Concerning the interaction between new mathematics and physics he said :
 
\begin{quote}

Mathematics and physics are more like Siamese twins. A separation in their case usually leads to death of at least one twin, so it should not be contemplated.
Some people sometimes oversimplify by regarding mathematics and physics as separate organisms, and by implication that mathematicians and physicists are different people. How does one decide which of the great figures of history was a mathematician or a physicist: Archimedes, Newton, Gauss, Hamilton, Maxwell, Riemann, Poincar\'e, Weyl?

Before Dirac's delta there was Heaviside (an engineer). Spinors are really due to Hamilton (both mathematician and physicist); \'Elie Cartan built on Sophus Lie; fibre bundles and U(1) gauge theory are due to Clifford. Grassmann used Maxwell's understanding of Hamilton's quaternions. Maxwell derived all the notation of div, grad and curl from Hamilton, who first wrote down the equation   
$${\left(i \frac{d}{dx} + j \frac{d}{dy}  +k \frac{d}{dz}\right)}^ 2  = - \Delta \,\,\,$$ and said that this must have deep physical meaning. This was many decades before Dirac. The dialogues continue and we now have a better understanding of fermions, supersymmetry and Morse theory. Also of topological quantum field theory in three-dimensions and JonesÕ knot invariants; Donaldson invariants in four dimensions, and the unique role of dimension four, both in physics and geometry. Also moduli spaces in Yang-Mills theory; monopoles and instantons; anomalies, cohomology and index theory; holography and geometry, etc. 
\end{quote}

\section{Understanding as a way of life}

\noindent In recent years, Atiyah concerned himself with relating our state of mind and  interest in mathematics to some cognitive notions like creativity, beauty, imagination and dreams \cite{kouneiher3, interview}.

\noindent For Atiyah creativity comes before any other act in mathematics, it is the triggering element. It is deeply embedded  in our mind.

\begin{quote}
People think mathematics begins when you write down a theorem followed by a proof. That's not the beginning, that's the end. For me the creative place in mathematics comes before you start to put things down on paper, before you try to write a formula. You picture various things, you turn them over in your mind. You are trying to create, just as a musician is trying to create music, or a poet. There are no rules laid down. You have to do it your own way. But at the end, just as a composer has to put it down on paper, you have to write things down. But the most important stage is understanding. A proof by itself doesn't give you understanding. You can have a long proof and no idea at the end of why it works. But to understand why it works, you have to have a kind of gut reaction to the thing. You have got to feel it.
\end{quote}

\noindent As for the beauty, as Dirac and Einstein and mostly Weyl (his hero along with Maxwell) Atiyah places the ``Beauty" at the heart of what we can call the essence of mathematics. It is not the kind of beauty that you can point to - it is beauty in a much more abstract sense.

\begin{quote}
It's been known for a long time that some part of the brain lights up when you listen to nice music, or read nice poetry, or look at nice pictures - and all of those reactions happen in the same place\footnote{ The Òemotional brain,Ó specifically the medial orbitofrontal cortex}. And the question was: Is the appreciation of mathematical beauty the same, or is it different? And the conclusion was, it is the same. The same bit of the brain that appreciates beauty in music, art and poetry is also involved in the appreciation of mathematical beauty. And that was a big discovery \cite{interview}.
\end{quote}

For Atiyah Euler' s equation $e^{i\pi} + 1 = 0$ embodies this beauty:

\begin{quote}
Euler's equation involves $\pi$; the mathematical constant Euler's number $e$ (2.71828..); $i$, the imaginary unit; $1$; and $0$.  It combines all the most important things in mathematics in one formula, and that formula is really quite deep. So everybody agreed that was the most beautiful equation. I used to say it was the mathematical equivalent of Hamlet's phrase ``To be, or not to be" -- very short, very succinct, but at the same time very deep. Euler's equation uses only five symbols, but it also encapsulates beautifully deep ideas, and brevity is an important part of beauty.
\end{quote}

\noindent Concerning the place of dreams in mathematics where ``the ideas appear fast and furious,  intuitive, imaginative, vague and clumsy commodities".

\begin{quote}
Dreams happen during the daytime, they happen at night. You can call them a vision or intuition. But basically they are a state of mind -- without words, pictures, formulas or statements. It is {\it pre} all that. It's pre-Plato. It is a very primordial feeling. And again, if you try to grasp it, it always dies. So when you wake up in the morning, some vague residue lingers, the ghost of an idea. You try to remember what it was and you only get half of it right, and maybe that is the best you can do \cite{interview}.
\end{quote}

\noindent In the D\'echiffreurs published in the Notices of the AMS, Atiyah wrote about the place of  dreams in our life:

\begin{quote}
In the broad light of day mathematicians check their equations and their proofs, leaving no
stone unturned in their search for rigour. But, at night, under the full moon, they dream, they
float among the stars and wonder at the miracle of the heavens. They are inspired. Without
dreams there is no art, no mathematics, no life\footnote{Atiyah M. F., {\it Les D\'echiffreurs} (2008), Notices of the AMS, (2010).}.
\end{quote}

\noindent Atiyah was a man in a hurry, his mind  crossed by an uninterrupted flow of teeming reflections and new ideas. He was a man who cogitates. Any exchange or discussion was a good opportunity to bring out mathematics and mathematical physical thinking. He enjoyed being around and  talking (even under the rain and for hours!) with colleagues and young researchers. He was a great conversationalist. 

 \noindent However, despite the tornado that could be released  by his presence and his willingness to communicate, interrogate and especially to be heard, Atiyah was paradoxically intrinsically a solitary man  who loved to be surrounded.  Atiyah was a man for whom friendship mattered. He did not have a superiority complex, but knew how to step back and advance people who  played a role in his mathematical life.

 \noindent  Atiyah fondly acknowledged the place of his wife Lily in his life and how she  allowed him the luxury of realizing his journey. She was a mathematician and left the academic path  so that he could go forward in his own career.

\noindent  But if there is one word that describes his life well, a constant feature of his story is ``to understand".  Atiyah was always guided by the desire to understand.

\begin{quote}
I always want to try to understand why things work. I'm not interested in getting a formula without knowing what it means. I always try to dig behind the scenes, so if I have a formula, I understand why it's there.
\end{quote}

\noindent Atiyah was always  a visionary and he was less interested in the past than in new ideas about the future.

\section{Acknowledgement}
We want to warmly thank Nigel Hitchin for his reading and suggestions during the preparation of the tribute, Erica Flapan for her support and Chika Mese for handling the tribute. We want to acknowledge the the refree's careful reading, remarks and suggestions.

\end{document}